\documentstyle[12pt]{amsart}
\newcommand{\g}{\goth}

\newcommand{\gtsl}{\mbox{\g sl}}

\newcommand{\hgtsl}{\mbox{$\hat{\gtsl}$}}

\newcommand{\nc}{\mbox{${\Bbb C}$}}
\newcommand{\nz}{\mbox{${\Bbb Z}$}}

\newcommand{\cF}{\mbox{${\cal F}$}}

\newcommand{\cP}{\mbox{${\cal P}$}}

\newcommand{\vep}{\varepsilon}
\newcommand{\vphi}{\varphi}
\theoremstyle{plain}
 \newtheorem{thm}{Theorem}[section]
 \newtheorem{prop}[thm]{Proposition}
 \newtheorem{lemma}[thm]{Lemma}
 
\theoremstyle{definition}
 \newtheorem{defn}{Definition}[section]

\theoremstyle{remark}
\newtheorem{rem}{Remark}

 \newtheorem{ack}{Acknowledgment}
      
\begin{document}
\title{Drinfeld comultiplication and vertex operators\\}
\author{Jintai Ding}
\author{kenji Iohara}
 
\address{Jintai Ding, RIMS, Kyoto University}

\address{Kenji Iohara, Dept. of Math., Kyoto University}

\maketitle
\begin{abstract}
For the current realization of the quantum affine algebras, Drinfeld gave a 
simple comultiplication of the quantum current operators. With this 
comultiplication, we study the related vertex operators for 
the case of $U_q(\hgtsl_n)$ and give an  
explicit bosonization of these new vertex operators. We use these 
vertex operators to construct the quantum current operators of 
 $U_q(\hgtsl_n)$ and discuss its connection with 
quantum boson-fermion correspondence. 
\end{abstract}
\pagestyle{plain}
\section{Introduction.} 

 Discovered by Drinfeld \cite{Dr1} 
and Jimbo \cite{J1}, 
Quantum group as a Hopf algebra  
presents  a new structure in both mathematics and physics. One of the 
fundamental features of this   structure comes from  its 
comultiplication, which is basically related  to all the new  concepts in 
this structure, such as R-matrix, etc. The  most  widely used
 comultiplication  is given  with  the definition of a  quantum group
 by the basic generators and the relations based on the data coming from the  
corresponding Cartan matrix. However, in  the case of quantum affine algebras, 
there is a  neglected aspect of the story. In this case, Drinfeld 
presented a different formulation of quantum affine algebras with generators 
in the form of current operators \cite{Dr3}, for which he proposed another 
comultiplication formula \cite{DF1} based  on the current
 formulation. The 
fundamental feature of this comultiplication is its simplicity, 
as opposed to the comultiplication formula induced from the conventional
 comultiplication which  
simply can not be written in a closed form with those current operators.
However, this new comultiplication is seldom  used beyond its definition. 
We propose to use this comultiplication formula to study 
the vertex operators coming from such a comultiplication.

In this paper,  we will first study the vertex operators  
for the fundamental 
representations at the level 1 for the case of $U_q(\hgtsl_n)$. 
Using  the  Frenkel-Jing construction of level 1 representations of 
 $U_q(\hgtsl_n)$ and the conventional comultiplication, the Kyoto 
group \cite{DFJMN} \cite{DO} \cite{Ko} studied
the related vertex operators and its application to 
XXZ model in statistical mechanics. They obtained
 the bosonization of the 
vertex operators, which  is partially  incomplete because    
the bosonization can be done only for one component of the vertex operators 
while the rest are implicit. With the new comultiplication, we 
give a complete and simple formula for the related intertwiners. In the 
classical case, these vertex operators are the corresponding Clifford 
algebras, from which we can reconstruct the level 1  representations, 
otherwise called spinor representations. This gives 
the  boson-fermion correspondence \cite{F1} \cite{Di1}. With the 
new formula of the vertex operators, which can be interpreted as the 
deformed fermions, we  recover the current operators of 
$U_q(\hgtsl_n)$. 

The main paper contains two sections. In Section 2, we will present the 
basic definitions and prove the new comultiplication 
formulas. In Section 3, we will present the Frenkel-Jing\cite{FJ}
 construction and 
present the explicit formulas for all the corresponding vertex operators. 
We will also give some formulas of the construction of the quantum 
current from those vertex operators and discuss its connection with 
quantum boson-fermion correspondence. In this paper, we will restrict 
ourselves to the case of $U_q(\hgtsl_n)$. The extension 
of these results to other cases is straightforward and will be the subject
in a subsequent paper. 

\section{Drinfeld comultiplication}
The first definition of quantum groups as 
a Hopf algebra given by Drinfeld \cite{Dr1} and Jimbo \cite{J1}  
are given in terms of 
basic generators and relations with corresponding Cartan matrices. 
For the case of quantum affine algebras, Drinfeld gave a realization 
in terms of operators in the form of current \cite{Dr3}. 
We will first present this  realization for the case of 
 $U_q(\hgtsl_n)$. 

Let $A=(a_{ij})$ be the Cartan matrix of type $A_{n-1}$.
\begin{defn}
The algebra $U_q(\hgtsl_n)$ is an associative algebra with unit 
1 and the generators: $\vphi_i(-m)$,$\psi_i(m)$, $x^{\pm}_i(l)$, for 
$i=i,...,n-1$, $l\in \nz $ and $m\in \nz_{\geq 0}$ and central
 elements $q^{\pm\frac{1}{2}c}$. Let $z$ be a formal variable and 
 $x_i^{\pm}(z)=\sum_{l\in \nz}x_i^{\pm}(l)z^{-l}$, 
$\vphi_i(z)=\sum_{m\in -\nz_{\geq 0}}\vphi_i(m)z^{-m}$ and 
$\psi_i(z)=\sum_{m\in \nz_{\geq 0}}\psi_i(m)z^{-m}$. In terms of the 
formal variables, 
the defining relations are 
\begin{align*}
& q^{\frac{1}{2}c}q{-\frac{1}{2}c}=1, \\
& \vphi_i(z)\vphi_j(w)=\vphi_j(w)\vphi_i(z), \\
& \psi_i(z)\psi_j(w)=\psi_j(w)\psi_i(z), \\
& \vphi_i(z)\psi_j(w)\vphi_i(z)^{-1}\psi_j(w)^{-1}=
  \frac{g_{ij}(z/wq^{-c})}{g_{ij}(z/wq^{c})}, \\
& \vphi_i(z)x_j^{\pm}(w)\vphi_i(z)^{-1}=
  g_{ij}(z/wq^{\mp \frac{1}{2}c})^{\pm1}x_j^{\pm}(w), \\
& \psi_i(z)x_j^{\pm}(w)\psi_i(z)^{-1}=
  g_{ij}(w/zq^{\mp \frac{1}{2}c})^{\mp1}x_j^{\pm}(w), \\
& [x_i^+(z),x_j^-(w)]=\frac{\delta_{i,j}}{q-q^{-1}}
  \left\{ \delta(z/wq^{-c})\psi_i(wq^{\frac{1}{2}c})-
          \delta(z/wq^{c})\vphi_i(zq^{\frac{1}{2}c}) \right\}, \\
& (z-q^{\pm a_{ij}}w)x_i^{\pm}(z)x_j^{\pm}(w)=
  (q^{\pm a_{ij}}z-w)x_j^{\pm}(w)x_i^{\pm}(z), \\
& [x_i^{\pm}(z),x_j^{\pm}(w)]=0 \quad \text{ for $a_{ij}=0$}, \\
& x_i^{\pm}(z_1)x_i^{\pm}(z_2)x_j^{\pm}(w)-(q+q^{-1})x_i^{\pm}(z_1)
  x_j^{\pm}(w)x_i^{\pm}(z_2)+x_j^{\pm}(w)x_i^{\pm}(z_1)x_i^{\pm}(z_2) \\
& +\{ z_1\leftrightarrow z_2\}=0, \quad \text{for $a_{ij}=-1$}
\end{align*}
where
\[ \delta(z)=\sum_{k\in \nz}z^k, \quad
   g_{ij}(z)=\frac{q^{a_{ij}}z-1}{z-q^{a_{ij}}}\quad \text{about $z=0$} \]
\end{defn}

In \cite{Dr3}, Drinfeld only gave the formulation of the algebra.
If we  extend the conventional comultiplication to 
those current operators,
the result would be a very  complicated formula which  can not be written 
in a closed form with only those current operators. 
However, Drinfeld also gave the Hopf algebra structure 
for such a formulation in an unpublished note \cite{DF1}.

\begin{thm}
The algebra $U_q(\hgtsl_n)$ has a Hopf algebra structure, which is  given 
by the following formulae. 

\noindent{\bf Coproduct $\Delta$}
\begin{align*}
\text{(0)}& \quad \Delta(q^{\frac c2})=q^{\frac c 2}\otimes q^{\frac c 2}, \\
\text{(1)}& \quad \Delta(x_i^+(z))=x_i^+(z)\otimes 1+
            \vphi_i(zq^{\frac{c_1}{2}})\otimes x_i^+(zq^{c_1}), \\
\text{(2)}& \quad \Delta(x_i^-(z))=1\otimes x_i^-(z)+
            x_i^-(zq^{c_2})\otimes \psi_i(zq^{\frac{c_2}{2}}), \\
\text{(3)}& \quad \Delta(\vphi_i(z))=
            \vphi_i(zq^{-\frac{c_2}{2}})\otimes\vphi_i(zq^{\frac{c_1}{2}}), \\
\text{(4)}& \quad \Delta(\psi_i(z))=
            \psi_i(zq^{\frac{c_2}{2}})\otimes\psi_i(zq^{-\frac{c_1}{2}}),
\end{align*}
where $q^{\frac {c_1}2}= q^{\frac c2}\otimes 1$ and 
 $q^{\frac {c_2}2}= 1\otimes q^{\frac c2}$. 
\noindent{\bf Counit $\vep$}
\begin{align*}
\vep(q^c)=1 & \quad \vep(\vphi_i(z))=\vep(\psi_i(z))=1, \\
            & \quad \vep(x_i^{\pm}(z))=0.
\end{align*}
\noindent{\bf Antipode $\quad a$}
\begin{align*}
\text{(0)}& \quad a(q^c)=q^{-c}, \\
\text{(1)}& \quad a(x_i^+(z))=-\vphi_i(zq^{-\frac{c}{2}})^{-1}
                               x_i^+(zq^{-c}), \\
\text{(2)}& \quad a(x_i^-(z))=-x_i^-(zq^{-c})
                               \psi_i(zq^{-\frac{c}{2}})^{-1}, \\
\text{(3)}& \quad a(\vphi_i(z))=\vphi_i(z)^{-1}, \\
\text{(4)}& \quad a(\psi_i(z))=\psi_i(z)^{-1}.
\end{align*}

\end{thm}
It is therefore 
clear that the comultiplication structure requires certain completion
on the tenor space. For certain  representations, such 
as the n-dimensional representations of $U_q(\hgtsl_n)$, which will be 
presented in the next section,
 this comultiplication may not be  well defined. Nevertheless, in the
case of 
any two highest weight representations,  this comultiplication
is well-defined and is already used in \cite{DM} to study the poles and zeros
of the current operators for integrable representations. 
We will further 
 present the  proof for the theorem above for the case of 
$U_q(\hgtsl_2)$.
However, this proof  should
be completely attributed to Drinfeld \cite{DF1}.

 \noindent\underline{Proof for the case of  $U_q(\hgtsl_2)$}

For the comultiplication above we have that 
\begin{align*}
\text{(1)}& \quad \Delta(x_1^+(z))=x_1^+(z)\otimes 1+
            \vphi_1(zq^{\frac{c_1}{2}})\otimes x_1^+(zq^{c_1}), \\
\text{(2)}& \quad \Delta(x_1^-(z))=1\otimes x_1^-(z)+
            x_1^-(zq^{c_2})\otimes \psi_1(zq^{\frac{c_2}{2}}), \\
\text{(3)}& \quad \Delta(\vphi_1(z))=
            \vphi_1(zq^{-\frac{c_2}{2}})\otimes\vphi_1(zq^{\frac{c_1}{2}}), \\
\text{(4)}& \quad \Delta(\psi_1(z))=
            \psi_1(zq^{\frac{c_2}{2}})\otimes\psi_1(zq^{-\frac{c_1}{2}}).
\end{align*}
It is clear that 
 $$\Delta\vphi_1(z)\Delta\vphi_1(w)=\Delta\vphi_1(w)\Delta\vphi_1(z),$$
$$\Delta\psi_1(z)\Delta\psi_1(w)=\Delta\psi_1(w)\Delta\psi_1(z).$$
Then, 
$$ \Delta\vphi_1(z)\Delta\psi_1(w)\Delta\vphi_1(z)^{-1}\Delta\psi_1(w)^{-1}=$$
  $$\vphi_1(zq^{-\frac{c_2}{2}})\psi_1(wq^{\frac{c_2}{2}})
\vphi_1(zq^{-\frac{c_2}{2}})^{-1}\psi_1(wq^{\frac{c_2}{2}})^{-1}
\otimes  \vphi_1(zq^{\frac{c_1}{2}})\psi_1(wq^{-\frac{c_1}{2}})
\vphi_1(zq^{\frac{c_1}{2}})^{-1}\psi_1(wq^{-\frac{c_1}{2}})^{-1}=$$
$$   \frac{g_{11}(z/wq^{-c_1}q^{-{c_2}} )}
{g_{11}(z/wq^{c_1}q^{-{c_2}})}
\frac{g_{11}(z/wq^{-c_2}q^{{c_1}})}{g_{11}(z/wq^{c_2}q^{{c_1}})}=$$
$$\frac{g_{11}(z/wq^{-c_1-c_2} )}
{g_{11}(z/wq^{c_1+c_2})}.$$

$$  \Delta\vphi_1(z)\Delta x_1^{+}(w)\Delta\vphi_1(z)^{-1}=$$
$$\vphi_1(zq^{-\frac{c_2}{2}})\otimes\vphi_1(zq^{\frac{c_1}{2}})
(x_1^+(w)\otimes 1+
            \vphi_1(wq^{\frac{c_1}{2}})\otimes x_1^+(wq^{c_1}))
 (\vphi_1(zq^{-\frac{c_2}{2}})\otimes\vphi_1(zq^{\frac{c_1}{2}}))^{-1}=$$
$$ g_{11}(z/wq^{ -\frac{1}{2}(c_1+c_2)})(x_1^+(w)\otimes 1+
            \vphi_1(wq^{\frac{c_1}{2}})\otimes x_1^+(wq^{c_1})). $$

$$  \Delta\vphi_1(z)\Delta x_1^{-}(w)\Delta\vphi_1(z)^{-1}=$$
$$\vphi_1(zq^{-\frac{c_2}{2}})\otimes\vphi_1(zq^{\frac{c_1}{2}})
(1\otimes x_1^-(w)+
            x_1^-(wq^{c_2})\otimes \psi_1(wq^{\frac{c_2}{2}}))
(\vphi_1(zq^{-\frac{c_2}{2}})\otimes\vphi_1(zq^{\frac{c_1}{2}}))^{-1}=$$
$$1\otimes x_1^-(w) g_{11}(z/wq^{ \frac{1}{2}(c_1+c_2)})^{-1}
+x_1^-(wq^{c_2})\otimes \psi_1(wq^{\frac{c_2}{2}})
g_{11}(z/wq^{ \frac{1}{2}(c_1-3c_2)})^{-1}
\frac {g_{11}(z/wq^{ \frac{1}{2}(c_1-3c_2)})}
      {g_{11}(z/wq^{ \frac{1}{2}(c_1+c_2)})}=$$
$$\Delta x_1^{-}(w)g_{11}(z/wq^{ \frac{1}{2}(c_1+c_2)})^{-1}.$$

The relation between $\Delta\psi_1(z)$ and $\Delta x_1^{\pm}(w)$
 can be proved in the same way demonstrated, 

$$ [\Delta x_1^+(z),\Delta x_1^-(w)]=$$
$$[x_1^+(z)\otimes 1+
            \vphi_1(zq^{\frac{c_1}{2}})\otimes x_1^+(zq^{c_1}),$$
$$ 
1\otimes x_1^-(w)+
            x_1^-(wq^{c_2})\otimes \psi_1(wq^{\frac{c_2}{2}})]=$$
$$[x_1^+(z)\otimes 1, x_1^-(wq^{c_2})\otimes \psi_1(wq^{\frac{c_2}{2}})]+
[\vphi_1(zq^{\frac{c_1}{2}})\otimes x_1^+(zq^{c_1}), 
1\otimes x_1^-(w)]+$$
$$ [\vphi_1(zq^{\frac{c_1}{2}})\otimes x_1^+(zq^{c_1}), 
            x_1^-(wq^{c_2})\otimes \psi_1(wq^{\frac{c_2}{2}})].$$
It is easy to show that the last term above is 0, therefore we have that
$$ [\Delta x_1^+(z),\Delta x_1^-(w)]=$$
$$\frac{(\delta(z/wq^{-c_1-c_2})\psi_1(wq^{\frac{c_1}{2}+c_2})-
\delta(z/wq^{c_1-c_2})\vphi_1(zq^{\frac{c_1}{2}}))\otimes 
\psi_1(wq^{\frac{c_2}{2}})}{q-q^{-1}}+$$
$$\frac{\vphi_1(zq^{\frac{c_1}{2}})\otimes
 ( \delta(z/wq^{c_1-c_2})\psi_1(wq^{\frac{1}{2}c_2})-
          \delta(z/wq^{c_1+c_2})\vphi_1(zq^{c_1+\frac{1}{2}c_2}))}{q-q^{-1}}=$$
$$ \frac{ \delta(z/wq^{-c_1-c_2})\psi_1(zq^{\frac{c_2}{2}+\frac{c_1+c_2}{2}})
\otimes\psi_1(zq^{-\frac{c_1}{2}+\frac{c_1+c_2}{2}})-
\delta(z/wq^{c_1+c_2})\vphi_1(zq^{-\frac{c_2}{2}+\frac{c_1+c_2}{2}})
\otimes\vphi_1(zq^{\frac{c_1}{2}+\frac{c_1+c_2}{2}})}{q-q^{-1}}.$$

$$ (z-q^{2}w)\Delta x_1^{+}(z)\Delta x_1^{+}(w)=$$
$$(z-q^2w)(x_1^+(z)\otimes 1+
            \vphi_1(zq^{\frac{c_1}{2}})\otimes x_1^+(zq^{c_1}))
(x_1^+(w)\otimes 1+
            \vphi_1(wq^{\frac{c_1}{2}})\otimes x_1^+(wq^{c_1}))=$$ 
$$ (q^{2}z-w)x_1 ^{+}(w)x_1^{+}(z)\otimes 1+
(z-q^2w)\vphi_1(wq^{\frac{c_1}{2}})x_1^+(z)\otimes x_1^+(wq^{c_1}) 
\left(\frac {q^2z/w-1}{z/w-q^2} \right)+$$
$$
(z-q^2w)x_1^+(w)\vphi_1(zq^{\frac{c_1}{2}})\otimes  
x_1^+(zq^{c_1}) \left(\frac {q^2z/w-1}{z/w-q^2}\right)+$$
$$ (q^{2}z-w) \vphi_1(wq^{\frac{c_1}{2}})\vphi_1(zq^{\frac{c_1}{2}})
   \otimes x_1^+(wq^{c_1})x_1^+(zq^{c_1})= $$
$$(q^{2}z-w)\Delta x_1^{+}(w)\Delta x_1^{+}(z).$$

Similarly we can prove the relation between 
$\Delta x^-(z)$ and $\Delta x^-(w)$.

Let $M$ be the operator from $U_q(\hgtsl_n)\otimes U_q(\hgtsl_n)$ to 
$U_q(\hgtsl_n)$ defined by the algebra multiplication. We can check that 
$M(1\otimes \vep)\Delta=\text{ id}$. 

$$ M(1\otimes a) \Delta(x_i^+(z))=$$
$$M(1\otimes a)( x_i^+(z)\otimes 1+
            \vphi_i(zq^{\frac{c_1}{2}})\otimes x_i^+(zq^{c_1}))=$$
$$ x_i^+(z)-\vphi_i(zq^{\frac{c}{2}})(\vphi_i(zq^{\frac{c}{2}}))^{-1}
 x_i^+(z)=$$
$$0=\vep(x_i^+(z))$$
Similarly we can check all the other relations to show that 
$M(a\otimes 1)=\vep.$ 
Thus, we prove that the comultiplication, the counit and the 
antipode give a Hopf algebra
structure for  $U_q(\hgtsl_2)$. 

Presently, there still exist a number of open problems
related to this new Hopf algebra structure. 
One is  whether this Hopf 
algebra is isomorphic, as a Hopf algebra,  to the 
conventional  $U_q(\hgtsl_n)$. 

\section{Bosonization of vertex operators}

\noindent{\bf vector representation}

Let $V=\oplus_{i=0}^{n-1} \nc|i\rangle$ be an n-dimensional space, 
where $\{ |i\rangle \}$ is its standard basis. Let 
$V_z=V\otimes \nc[z,z^{-1}]$, where $z$ is a formal variable.

\begin{lemma}There exists an  n-dimensional representation of $U_q(\hgtsl_n)$
on $V_z$. The action of the current operators is given by the following:
\begin{align*}
\circ& \quad x_i^+(w).|j\rangle= \delta_{ij}
                                 \delta(\frac{w}{q^iz})|i-1\rangle, \\
\circ& \quad x_i^-(w).|j\rangle= \delta_{i-1j}
                                 \delta(\frac{w}{q^iz})|i\rangle, \\
\circ& \quad \begin{cases}
             \vphi_i(w).|i-1\rangle 
       & =\dfrac{q^{-1}-q^{-i+1}\frac{w}{z}}{1-q^{-i}\frac{w}{z}}
         |i-1\rangle, \\
             \vphi_i(w).|i\rangle
       & =\dfrac{q-q^{-i-1}\frac{w}{z}}{1-q^{-i}\frac{w}{z}}
         |i\rangle    \\
             \vphi_i(w).|j\rangle
       &=|j\rangle \quad (j\neq i,i-1)
       \end{cases} \in \nc[[\frac{w}{z}]], \\
\circ& \quad \begin{cases}
             \psi_i(w).|i-1\rangle
       & =\dfrac{q-q^{i-1}\frac{z}{w}}{1-q^{i}\frac{z}{w}}
         |i-1\rangle, \\
             \psi_i(w).|i\rangle
       & =\dfrac{q^{-1}-q^{i+1}\frac{z}{w}}{1-q^{i}\frac{z}{w}}
         |i\rangle   \\
             \psi_i(w).|j\rangle
       & =|j\rangle \quad (j\neq i,i-1)
       \end{cases} \in \nc[[\frac{z}{w}]].
\end{align*}
where we set $|-1\rangle =|n\rangle=0$.
\end{lemma} 

Let $E_{ij}~(i,j=0,1)$ be the matrix elements such that 
$E_{ij}|k\rangle =\delta_{kj}|i\rangle$, let $R(z/w)$ be an 
operator on $V_z\otimes V_w$, which is defined as: 
$$E_{00}\otimes E_{00}+E_{11}\otimes E_{11}+ E_{00}\otimes E_{11}
\frac{q-q^{-1}z/w}{1-z/w}+E_{11}\otimes E_{00} \frac{1-z/w} {q^{-1}-qz/w}.$$

\begin{prop}
Let $x$ be an operator in $U_q(\hgtsl_2)$, then  on $V_z\otimes V_w$, we have
that $R(z/w)\Delta(x)=\Delta'(x)R(z/w),$ 
where $\Delta'$ is the opposite comultiplication. 
\end{prop}

Similarly, we can write the operator $R(z/w)$ for $U_q(\hgtsl_n)$, which is 
diagonal. 

Next, we will give the Frenkel-Jing construction of level $1$ representation 
of  $U_q(\hgtsl_n)$ on the Fock space. 

Let $\overline{Q}=\oplus_{j=1}^{n-1} \nz \alpha_j$ be the root lattice of 
$\gtsl_n$, $\overline{\Lambda}_j=\Lambda_j-\Lambda_0$ be the classical 
part of the $i$-th fundamental weight.

Let Heisenberg algebra be an algebra generated by 
$\{ a_{i,k}|1\leq i \leq n-1,~k\in \nz \setminus \{ 0\} \}$ 
satisfying:
\[ [a_{i,k},a_{j,l}]=\delta_{k+l,0}\frac{[(\alpha_i,\alpha_j)k][k]}{k}. \]

Now let's define a group algebra $\nc(q) [\overline{\cP}]$. 
Let $\overline{\cP}$ be the weight lattice of $\gtsl_n$. 
We fix our free basis $\alpha_2,\cdots, \alpha_{n-1},\overline{\Lambda}_{n-1}$.
They satisfy
\begin{align*}
& e^{\alpha_i}e^{\alpha_j}=(-1)^{(\alpha_i,\alpha_j)}e^{\alpha_j}e^{\alpha_i}
  \qquad 2\leq i,j\leq n-1, \\
& e^{\alpha_i}e^{\overline{\Lambda}_{n-1}}=(-1)^{\delta_{i,n-1}}
  e^{\overline{\Lambda}_{n-1}}e^{\alpha_i}, 
  \qquad 2\leq i\leq n-1.
\end{align*}
For $\alpha=m_2\alpha_2+\cdots m_{n-1}\alpha_{n-1}
     +m_n\overline{\Lambda}_{n-1}$, we set
\[ e^{\alpha}=e^{m_2\alpha_2}\cdots e^{m_{n-1}\alpha_{n-1}}
              e^{m_n\overline{\Lambda}_{n-1}}. \]
Note that the following equations hold.
\begin{align*}
& \overline{\Lambda}_i=-\alpha_{i+1}-2\alpha_{i+2}-\cdots 
                      -(n-i-1)\alpha_{n-1}+(n-i)\overline{\Lambda}_{n-1}, \\
& \alpha_1=-2\alpha_2-3\alpha_3-\cdots -(n-1)\alpha_{n-1}
           +n\overline{\Lambda}_{n-1}.
\end{align*}

Put
\[ a_{ik}^{\ast}=\frac{1}{[k]^2[nk]}\sum_{j=1}^{n-1}
   [ \min (i,j)k][ \min (n-i,n-j)k]a_{jk},\quad 1\leq i< n, k\neq 0. \]
Then they satisfy
\[ [ a_{ik}^{\ast}, a_{jl}]= \delta_{ij}\delta_{k+l,0}\frac{[k]}{k}. \]
We define the Fock space as:
\[ \cF_{i}:=\nc(q) [a_{j,-k}(1\leq j\leq n-1,~k\in \nz_{>0})]\otimes
            \nc(q) [\overline{Q}]
              e^{\overline{\Lambda}_i}
   \quad (0\leq i\leq n-1). \]

The action of operators 
$a_{j,k},\partial_{\alpha_j},e^{\alpha}$~
$(1\leq j \leq n-1, \alpha \in \overline{Q})$ is given by 
\begin{align*}
a_{j,k}\cdot f\otimes e^{\beta}& =\begin{cases}
                      a_{j,k}f \otimes e^{\beta}          & k< 0  \\
            \text{$[a_{j,k},f]$}  \otimes e^{\beta} \quad & k> 0
                                  \end{cases}, \\
\partial_{\alpha}\cdot f\otimes e^{\beta}& 
        =(\alpha,\beta)f\otimes e^{\beta}, 
         \qquad \text{for}~f\otimes e^{\beta}\in \cF_{i} \\
e^{\alpha}\cdot f\otimes e^{\beta} & 
=f\otimes e^{\alpha}e^{\beta}.
\end{align*}

\begin{lemma} 
Let 
\begin{align*}
\circ & \quad x_j^{\pm}(z)\mapsto
        \exp[\pm \sum_{k>0}\frac{a_{j,-k}}{[k]}q^{\mp \frac{1}{2}k}z^k]
        \exp[\mp \sum_{k>0}\frac{a_{j,k}}{[k]}q^{\mp \frac{1}{2}k}z^{-k}]
        e^{\pm \alpha_j}z^{\pm \partial_{\alpha_j}+1}, \\
\circ & \quad \vphi_j(z)\mapsto
        \exp[-(q-q^{-1})\sum_{k>0}a_{j,-k}z^k]q^{-\partial_{\alpha_j}}, \\
\circ & \quad \psi_j(z)\mapsto
        \exp[(q-q^{-1})\sum_{k>0}a_{j,k}z^{-k}]q^{\partial_{\alpha_j}},
\end{align*}
which  gives level 1 highest representations with  the $i$-th  fundamental
weight on $\cF_{i}$. 
\end{lemma}

\begin{defn}
Vertex operators are intertwiners of the following types: 
\begin{enumerate}
\renewcommand{\labelenumi}{(\roman{enumi})}
\item (type $I$) \hspace{1 in} 
      $\displaystyle{ \Phi^{(i,i+1)}(z):\cF_{i+1}\longrightarrow
                      \cF_{i}\otimes V_{z}, }$ 
\item (type $II$) \hspace{0.95 in}
      $\displaystyle{ \Psi^{(i,i+1)}(z):\cF_{i+1}\longrightarrow
                      V_{z}\otimes \cF_{i}. }$
\item (dual of type $I$) \hspace{0.5 in}
      $\displaystyle{ \Phi^{*(i+1,i)}(z):\cF_{i}\otimes V_{z}
                      \longrightarrow \cF_{i+1}, }$
\item (dual of type $II$)\hspace{0.45 in}
      $\displaystyle{ \Psi^{*(i+1,i)}(z):V_{z}\otimes \cF_{i}
                      \longrightarrow \cF_{i+1}, }$
\end{enumerate}
Here the indices are considered modulo $n$.
\end{defn}
Set
\begin{align*}
&  \Phi^{(i,i+1)}(z)=\sum_{j=0}^{n-1} \Phi^{(i.i+1)}_j(z)\otimes |j\rangle,
\qquad
   \Psi^{(i,i+1)}(z)=\sum_{j=0}^{n-1} |j\rangle\otimes \Psi^{(i.i+1)}_j(z),
\\
&  \Phi^{*(i+1,i)}(z)(u\otimes |j\rangle )=\Phi^{*(i+1,i)}_j(z)u, \quad
   \Psi^{*(i+1,i)}(z)(|j\rangle \otimes u)=\Psi^{*(i+1,i)}_j(z)u, \qquad
   u \in \cF_{i+1}.
\end{align*}
The normalization of those operators is  given by 
\begin{enumerate}
\renewcommand{\labelenumi}{(\roman{enumi})}
\item $\displaystyle{ 
       \langle \Lambda_{i}|\Phi^{(i,i+1)}_{i}(z)|\Lambda_{i+1}\rangle=1,
\qquad \langle \Lambda_{i}|\Psi^{(i,i+1)}_{i}(z)|\Lambda_{i+1}\rangle=1,    
       } $
\item $\displaystyle{
       \langle \Lambda_{i+1}|\Phi^{*(i+1,i)}_{i}(z)|\Lambda_{i}\rangle=1,
\qquad \langle \Lambda_{i+1}|\Psi^{*(i+1,i)}_{i}(z)|\Lambda_{i}\rangle=1.
       } $
\end{enumerate}

We will first give the list of OPEs (operator product expansions), where 
we abbreviate the superscript of vertex operators.

\noindent \underline{Type $I,II$}
\begin{align*}
& [\Phi_j(z),x_i^+(w)]=\delta_{i-1,j}\delta(\frac{w}{q^{i-1}z})
                       \vphi_i(wq^{\frac{1}{2}})\Phi_i(z), \\
& \Phi_j(z)x_i^-(w)=\begin{cases}
                    \dfrac{q-q^{i-1}\frac{z}{w}}{1-q^i\frac{z}{w}}
                    x_i^-(w)\Phi_{i-1}(z) 
                    & \quad \text{for $j=i-1$} \\
                    \dfrac{q^{-1}-q^{i+1}\frac{z}{w}}{1-q^i\frac{z}{w}}
                    x_i^-(w)\Phi_{i}(z)+\delta(\frac{w}{q^{i}z})
                    \Phi_{i-1}(z)
                    & \quad \text{for $j=i$} \\
                    x_i^-(w)\Phi_j(z)
                    & \quad \text{otherwise}
                    \end{cases}, \\
& \Phi_j(z)\vphi_i(w)=\begin{cases}
                    \dfrac{q^{-1}-q^{-i+\frac{3}{2}}\frac{w}{z}}
                          {1-q^{-i+\frac{1}{2}}\frac{w}{z}}
                    \vphi_i(w)\Phi_{i-1}(z)
                    & \quad \text{for $j=i-1$} \\
                    \dfrac{q-q^{-i-\frac{1}{2}}\frac{w}{z}}
                          {1-q^{-i+\frac{1}{2}}\frac{w}{z}}
                    \vphi_i(w)\Phi_i(z)
                    & \quad \text{for $j=i$} \\
                    \vphi_i(w)\Phi_j(z)
                    & \quad \text{otherwise}
                    \end{cases}, \\
& \Phi_j(z)\psi_i(w)=\begin{cases}
                    \dfrac{q-q^{i-\frac{1}{2}}\frac{z}{w}}
                          {1-q^{i+\frac{1}{2}}\frac{z}{w}}
                    \psi_i(w)\Phi_{i-1}(z) 
                    & \quad \text{for $j=i-1$} \\
                    \dfrac{q^{-1}-q^{i+\frac{3}{2}}\frac{z}{w}}
                          {1-q^{i+\frac{1}{2}}\frac{z}{w}}
                    \psi_i(w)\Phi_i(z)
                    & \quad \text{for $j=i$} \\
                    \psi_i(w)\Phi_j(z)
                    & \quad \text{otherwise}
                    \end{cases}, \\
& \Psi_j(z)x_i^+(w)=\begin{cases}
                    \dfrac{q^{-1}-q^{-i+1}\frac{w}{z}}{1-q^{-i}\frac{w}{z}}
                    x_i^+(w)\Psi_{i-1}(z)+\delta(\frac{w}{q^iz})
                    \Psi_i(z)
                    & \quad \text{for $j=i-1$} \\
                    \dfrac{q-q^{-i-1}\frac{w}{z}}{1-q^{-i}\frac{w}{z}}
                    x_i^+(w)\Psi_i(z)
                    & \quad \text{for $j=i$} \\
                    x_i^+(w)\Psi_j(z)
                    & \quad \text{otherwise}
                    \end{cases}, \\
& [\Psi_j(z),x_i^-(w)]=\delta_{i,j}\delta(\frac{w}{q^{i-1}z})
                       \psi_i(wq^{\frac{1}{2}})\Psi_{i-1}(z), \\
& \Psi_j(z)\vphi_i(w)=\begin{cases}
                    \dfrac{q^{-1}-q^{-i+\frac{1}{2}}\frac{w}{z}}
                          {1-q^{-i-\frac{1}{2}}\frac{w}{z}}
                    \vphi_i(w)\Psi_{i-1}(z)
                    & \quad \text{for $j=i-1$} \\
                    \dfrac{q-q^{-i-\frac{3}{2}}\frac{w}{z}}
                          {1-q^{-i-\frac{1}{2}}\frac{w}{z}}
                    \vphi_i(w)\Psi_i(z)
                    & \quad \text{for $j=i$} \\
                    \vphi_i(w)\Psi_j(z)
                    & \quad \text{otherwise}
                    \end{cases}, \\
& \Psi_j(z)\psi_i(w)=\begin{cases}
                    \dfrac{q-q^{i-\frac{3}{2}}\frac{z}{w}}
                          {1-q^{i-\frac{1}{2}}\frac{z}{w}}
                    \psi_i(w)\Psi_{i-1}(z)
                    & \quad \text{for $j=i-1$} \\
                    \dfrac{q^{-1}-q^{i+\frac{1}{2}}\frac{z}{w}}
                          {1-q^{i-\frac{1}{2}}\frac{z}{w}}
                    \psi_i(w)\Psi_i(z)
                    & \quad \text{for $j=i$} \\
                    \psi_i(w)\Psi_j(z)
                    & \quad \text{otherwise}
                    \end{cases},
\end{align*}
\noindent \underline{Dual of type $I,II$}
\begin{align*}
& [x_i^+(w),\Phi^*_j(z)]=\delta_{i,j}\delta(\frac{w}{q^{i-1}z})
                         \Phi^*_{i-1}(z)\vphi_i(wq^{\frac{1}{2}}), \\
& x_i^-(w)\Phi^*_j(z)=\begin{cases}
                     \dfrac{q-q^{i-1}\frac{z}{w}}{1-q^i\frac{z}{w}}
                     \Phi^*_{i-1}(z)x_i^-(w)+\delta(\frac{w}{q^iz})
                     \Phi^*_i(z)
                     & \quad \text{for $j=i-1$} \\
                     \dfrac{q^{-1}-q^{i+1}\frac{z}{w}}{1-q^i\frac{z}{w}}
                     \Phi^*_i(z)x_i^-(w)
                     & \quad \text{for $j=i$} \\
                     \Phi^*_j(z)x_i^-(w)
                     & \quad \text{otherwise}
                     \end{cases}, \\
& \vphi_i(w)\Phi^*_j(z)=\begin{cases}
                     \dfrac{q^{-1}-q^{-i+\frac{3}{2}}\frac{w}{z}}
                           {1-q^{-i+\frac{1}{2}}\frac{w}{z}}
                     \Phi^*_{i-1}(z)\vphi_i(w)
                     & \quad \text{for $j=i-1$} \\
                     \dfrac{q-q^{-i-\frac{1}{2}}\frac{w}{z}}
                           {1-q^{-i+\frac{1}{2}}\frac{w}{z}}
                     \Phi^*_i(z)\vphi_i(w)
                     & \quad \text{for $j=i$} \\
                     \Phi^*_j(z)\vphi_i(w)
                     & \quad \text{otherwise}
                     \end{cases}, \\
& \psi_i(w)\Phi^*_j(z)=\begin{cases}
                     \dfrac{q-q^{i-\frac{1}{2}}\frac{z}{w}}
                           {1-q^{i+\frac{1}{2}}\frac{z}{w}}
                     \Phi^*_{i-1}(z)\psi_i(w)
                     & \quad \text{for $j=i-1$} \\
                     \dfrac{q^{-1}-q^{i+\frac{3}{2}}\frac{z}{w}}
                           {1-q^{i+\frac{1}{2}}\frac{z}{w}}
                     \Phi^*_i(z)\psi_i(w)
                     & \quad \text{for $j=i$} \\
                     \Phi^*_j(z)\psi_i(w)
                     & \quad \text{otherwise}
                     \end{cases}, \\
& x_i^+(w)\Psi^*_j(z)=\begin{cases}
                     \dfrac{q^{-1}-q^{-i+1}\frac{w}{z}}{1-q^{-i}\frac{w}{z}}
                     \Psi^*_{i-1}(z)x_i^+(w)
                     & \quad \text{for $j=i-1$} \\
                     \dfrac{q-q^{-i-1}\frac{w}{z}}{1-q^{-i}\frac{w}{z}}
                     \Psi^*_i(z)x_i^+(w)+\delta(\frac{w}{q^iz})
                     \Psi^*_{i-1}(z)
                     & \quad \text{for $j=i$} \\
                     \Psi^*_j(z)x_i^+(w)
                     & \quad \text{otherwise}
                     \end{cases}, \\
& [x_i^-(w),\Psi^*_j(z)]=\delta_{i-1,j}\delta(\frac{w}{q^{i-1}z})
                        \Psi^*_i(z)\psi_i(wq^{\frac{1}{2}}), \\
& \vphi_i(w)\Psi^*_j(z)=\begin{cases}
                     \dfrac{q^{-1}-q^{-i+\frac{1}{2}}\frac{w}{z}}
                           {1-q^{-i-\frac{1}{2}}\frac{w}{z}}
                     \Psi^*_{i-1}(z)\vphi_i(w)
                     & \quad \text{for $j=i-1$} \\
                     \dfrac{q-q^{-i-\frac{3}{2}}\frac{w}{z}}
                           {1-q^{-i-\frac{1}{2}}\frac{w}{z}}
                     \Psi^*_i(z)\vphi_i(w)
                     & \quad \text{for $j=i$} \\
                     \Psi^*_j(z)\vphi_i(w)
                     & \quad \text{otherwise}
                     \end{cases}, \\
& \psi_i(w)\Psi^*_j(z)=\begin{cases}
                     \dfrac{q-q^{i-\frac{3}{2}}\frac{z}{w}}
                           {1-q^{i-\frac{1}{2}}\frac{z}{w}}
                     \Psi^*_{i-1}(z)\psi_i(w)
                     & \quad \text{for $j=i-1$} \\
                     \dfrac{q^{-1}-q^{i+\frac{1}{2}}\frac{z}{w}}
                           {1-q^{i-\frac{1}{2}}\frac{z}{w}}
                     \Psi^*_i(z)\psi_i(w)
                     & \quad \text{for $j=i$} \\
                     \Psi^*_j(z)\psi_i(w)
                     & \quad \text{otherwise}
                     \end{cases}.
\end{align*}

\begin{rem}We remark that 2 OPE's related to 
$\Psi_i(z)x_i^-(w)$ and $x_i^+(w)\Phi^*_i(z)$ are 
$$\Psi_i(z)x_i^-(w)=x_i^-(w)\Psi_i(z),$$
$$x_i^+(w)\Phi^*_i(z)=\Phi^*_i(z)x_i^+(w).$$
However the two sides of the equations are in different expansion directions, 
which implies that both sides have poles . 
These two equations will not degenerate the classical 
intertwiner relations when 
$q=1$. For this 
case, we impose the following locality conditions:
\[ (1-\frac{w}{q^{i+1}z})[\Psi_i(z),x_i^-(w)]=0, \qquad
   (1-\frac{w}{q^{i+1}z})[x_i^+(w),\Phi^*_i(z)]=0,  \]
which ensure that the degeneration of these  formulae is consistent with its 
classical degeneration when $q=1$. 
\end{rem}
 
\begin{thm}
The vertex operators are given as:
\begin{align*}
\Phi_j^{(i,i+1)}(z)=
& \exp[-\sum_{k>0}a_{j-k}^{\ast}q^{\frac{3}{2}k}(q^jz)^k
       +\sum_{k>0}a_{j+1-k}^{\ast}q^{\frac{1}{2}k}(q^jz)^k] \\
& \exp[-\sum_{k>0}a_{jk}^{\ast}q^{-\frac{1}{2}k}(q^jz)^{-k}
       +\sum_{k>0}a_{j+1k}^{\ast}q^{\frac{1}{2}k}(q^jz)^{-k}] \\
& e^{\overline{\Lambda}_j-\overline{\Lambda}_{j+1}}
  (q^{j+1}z)^{\partial_{\overline{\Lambda}_j}}
  (q^jz)^{-\partial_{\overline{\Lambda}_{j+1}}}
  (-1)^{(n-1)\partial_{\overline{\Lambda}_1}}(c)_j^i, \\
\Psi_j^{(i,i+1)}(z)=
& \exp[-\sum_{k>0}a_{j-k}^{\ast}q^{\frac{1}{2}k}(q^jz)^k
       +\sum_{k>0}a_{j+1-k}^{\ast}q^{-\frac{1}{2}k}(q^jz)^k] \\
& \exp[-\sum_{k>0}a_{jk}^{\ast}q^{-\frac{3}{2}k}(q^jz)^{-k}
       +\sum_{k>0}a_{j+1k}^{\ast}q^{-\frac{1}{2}k}(q^jz)^{-k}] \\
& e^{\overline{\Lambda}_j-\overline{\Lambda}_{j+1}}
  (q^{j+1}z)^{\partial_{\overline{\Lambda}_j}}
  (q^jz)^{-\partial_{\overline{\Lambda}_{j+1}}}
  (-1)^{(n-1)\partial_{\overline{\Lambda}_1}}(c)_j^i, \\
\Phi_j^{*(i+1,i)}(z)=
& \exp[\sum_{k>0}a_{j-k}^{\ast}q^{\frac{3}{2}k}(q^jz)^k
       -\sum_{k>0}a_{j+1-k}^{\ast}q^{\frac{1}{2}k}(q^jz)^k] \\
& \exp[\sum_{k>0}a_{jk}^{\ast}q^{-\frac{1}{2}k}(q^jz)^{-k}
       -\sum_{k>0}a_{j+1k}^{\ast}q^{\frac{1}{2}k}(q^jz)^{-k}] \\
& e^{-\overline{\Lambda}_j+\overline{\Lambda}_{j+1}}
  (q^{j+1}z)^{-\partial_{\overline{\Lambda}_j}}
  (q^jz)^{\partial_{\overline{\Lambda}_{j+1}}}
  (-1)^{-(n-1)\partial_{\overline{\Lambda}_1}}(c^{\ast})_j^i, \\
\Psi_j^{*(i+1,i)}(z)=
& \exp[\sum_{k>0}a_{j-k}^{\ast}q^{\frac{1}{2}k}(q^jz)^k
       -\sum_{k>0}a_{j+1-k}^{\ast}q^{-\frac{1}{2}k}(q^jz)^k] \\
& \exp[\sum_{k>0}a_{jk}^{\ast}q^{-\frac{3}{2}k}(q^jz)^{-k}
       -\sum_{k>0}a_{j+1k}^{\ast}q^{-\frac{1}{2}k}(q^jz)^{-k}] \\
& e^{-\overline{\Lambda}_j+\overline{\Lambda}_{j+1}}
  (q^{j+1}z)^{-\partial_{\overline{\Lambda}_j}}
  (q^jz)^{\partial_{\overline{\Lambda}_{j+1}}}
  (-1)^{-(n-1)\partial_{\overline{\Lambda}_1}}(c^{\ast})_j^i, 
\end{align*}
where $(c)_j^i=(-q)^{j-i}(c)_i^i,(c^{\ast})_j^i=(c^{\ast})_i^i$ and
\begin{align*}
(c)_i^i=&
[(-1)^{-(n-1)}z]^{\frac{n-i-1}{n}}(-1)^{\frac{1}{2}(n-i-1)(n-i-2)}, \\
(c^{\ast})_i^i=&
[(-1)^{n-1}]^{\frac{n-i}{n}}[q^nz]^{\frac{i}{n}}(-1)^{\frac{1}{2}(n-i)(n-i-1)}.
\end{align*}
Here $\ast$ means dual.
\end{thm} 

From this theorem, we can derive all the correlation functions of 
vertex operators. We will give a simple example of the case of 
$U_q(\hgtsl_2)$ to show what they are like. 

Let $v$ be the highest weigh vector in $\cF_{0}$
$$<v, \Phi_0^{(0,1)}(z_1) \Phi_0^{(1,0)}(z_2), v>= 
  <v, \Phi_1^{(0,1)}(z_1) \Phi_1^{(1,0)}(z_2), v>= 0,$$ 
$$<v, \Phi_1^{(0,1)}(z_2) \Phi_0^{(1,0)}(z_1), v>=-q^{-1}(z_1/z_2)^{1/2}
\frac {(qz_1/z_2;q^4)_\infty}{(q^3z_1/z_2,q^4)_\infty},$$
$$<v, \Phi_0^{(0,1)}(z_2) \Phi_1^{(1,0)}(z_1), v>=
z_1^{1/2}z_2^{3/2}\frac {(q^4z_1/z_2;q^4)_\infty}{(q^6z_1/z_2,q^4)_\infty},$$
where  $(z,p)_\infty=\prod^\infty_0(1-zp^m)$. 

With those vertex operators, we can derive the following 
formulae. 
\begin{thm}
\begin{align*}
\text{$(1)$}& \quad  :\Phi_{i-1}^{(i-1,i)}(z)\Phi_i^{\ast(i+1,i)}(z):
               =x_i^-(q^iz), \\ 
\text{$(2)$}& \quad  :\Psi_i^{(i-1,i)}(z)\Psi_{i-1}^{\ast(i+1,i)}(z):
               =-qx_i^+(q^iz), \\
\text{$(3)$}& \quad  :x_i^+(q^{\frac{1}{2}}z)x_i^-(q^{-\frac{1}{2}}z):
               =z^2 \psi_i(z), \\
\text{$(4)$}& \quad  :x_i^+(q^{-\frac{1}{2}}z)x_i^-(q^{\frac{1}{2}}z):
               =z^2 \vphi_i(z),
\end{align*}
where $: \cdot :$ is the normal ordering.
\end{thm}

In \cite{Di1}\cite{Di2}, we showed  that the vertex operators in the classical 
case have the structure of Clifford algebras, which is used to 
construct representations of $\hgtsl_n$. This is the 
boson-fermion correspondence. We also established  
a quantum version of boson-fermion correspondence, where the quantum 
fermions are identified as vertex operators with the conventional 
comultiplication. However the quantum boson-fermion correspondence is very 
implicit due to the incomplete formulae for intertwiners. With Theorem 
3.4 and Theorem 3.5, the correlation functions of all the 
vertex operators are clear. Utilizing 
 all these results,  we hope that we can 
present a different definition of quantum Clifford algebras and derive a 
more explicit boson-fermion correspondence. Furthermore
 we recognize that 
our formulae may be related to the quantum ${\cal W}$-algebras \cite{AKOS}
\cite{FF} especially the bosonization of those algebras.  

\begin{ack}
J.D. is supported by the grant Reward research (A) 08740020 from the 
Ministry of Education of Japan. K.I. is partly supported by the Japan Society
for the Promotion of Science. 
\end{ack}


\begin{thebibliography}{[Beck]}


\bibitem
  [AKOS] {AKOS}H. Awata, H. Kubo, S. Odake, J. Shiraishi
{\it Quantum ${\cal W}_N$ algebras and Macdonald Polynomial},
q-alg/9508011


\bibitem
 [DO] {DO}E. Date and M. Okado 
{\it Calculation of excitation spectra
of the spin model related with the vector representation of the
quantized affine algebra of type $A^{(1)}_n$},
 Int. J. of Mod. Phys.A, {\bf 9}, No.3, 1994


\bibitem     [DFJMN]
{DFJMN} B. Davis, O. Foda, M. Jimbo, T. Miwa and A. Nakayashiki
{\it Diagonalization of the XXZ Hamiltonian by
vertex operators}, CMP, 
{\bf 151}, 1993, 89-153  

\bibitem   [Di1] {Di1} J.Ding
{\it Spinor representations of $U_q(\hat{\frak gl}(n))$ 
and quantum Boson-Fermion correspondence}, 
RIMS-1043, q-alg/9510014


\bibitem   [Di2] {Di2} J.Ding
{\it  Spinor Representations of 
$U_q(\hat {\frak o}(N))$},
To appear in LMP, RIMS-1060, q-alg/9602021




\bibitem   [DF1] {DF1} J. Ding, I. B. Frenkel
{\it Isomorphism of two realizations of quantum affine algebra 
$U_q(\hat {\frak gl}(n))$ },
 CMP, {\bf 156}, 1993, 277-300
Physics 


\bibitem    [DF2] {DF2} J. Ding and I. B. Frenkel
 {\it Spinor and oscillator
representations of quantum groups, in:  Lie Theory and Geometry
in Honor of Bertram Kostant, Progress in mathematics}, 
{\bf 123}, Birkhauser, Boston,  1994

\bibitem [DM] {DM} J.Ding and T. Miwa 
{\it Zeros and poles of quantum current operators and 
the quantum integrable condition}, RIMS-1092. 


\bibitem    [Dr1] {Dr1} V. G. Drinfeld 
{\it Hopf algebra and the quantum Yang-Baxter Equation},
Dokl. Akad. Nauk. SSSR, {\bf 283}, 1985, 1060-1064

 

\bibitem    [Dr2] {Dr2} V.G. Drinfeld
{\it Quantum Groups},
ICM Proceedings, New York, Berkeley, 1986, 798-820

 

\bibitem  [Dr3]{Dr3}  V. G. Drinfeld {\it
 New realization of Yangian and quantum
affine algebra}, Soviet Math. Doklady,
{\bf 36}, 1988,  212-216


\bibitem 
 [FF] {FF}B. Feigin, E. Frenkel
{\it Quantum ${\cal W}$-algebra and elliptic algebras},
RIMS-1027
 

\bibitem 
 [F1] {F1} I. B. Frenkel 
{\it Two constructions of affine Lie algebra representations and
boson-fermion correspondence in quantum field theory},
 J. Funct. Anal.,{\bf  44},1981,259-327
 

\bibitem    [FJ]{FJ} I. B. Frenkel, N. Jing
{\it Vertex representations of quantum affine algebras},
 Proc. Natl. Acad. Sci., USA,{\bf 85}, 1988, 9373-9377



\bibitem     [J1]{J1} M. Jimbo 
{\it A $q$-difference analogue of $U({\frak g})$ and Yang-Baxter equation
},
Lett. Math. Phys. {\bf 10}, 1985, 63-69



\bibitem     [Ko]{Ko} Y. Koyama
{\it Staggered Polarization of Vertex Model
with $U_q(\hat{sl}(n))$-symmetry
},
CMP, {\bf 164}, 1994, 277-291

\end{thebibliography}
\end{document}